\documentclass[12pt]{article}
\usepackage{amsmath,amsfonts,amssymb,amsthm}

\usepackage{fullpage}
\def\half{\frac{1}{2}}
\newcommand{\7}{\dagger}            
\newcommand{\up}{{\mathord{\uparrow}}} 
\newcommand{\dn}{{\mathord{\downarrow}}} 
\newcommand{\ox}{\otimes}           
\newcommand{\lt}{\triangleright}    
\newcommand{\rt}{\triangleleft}     

\def\H{{\cal H}}
\def\Z{\Bbb Z}
\def\N{\Bbb N}

\def\IC{\Bbb C}
\def\II{\Bbb I}
\def\e{\varepsilon}
\def\ee{\epsilon}
\def\a{\alpha}

\def\b{\beta}
\def\g{\gamma}
\def\d{\delta}

\newcommand{\sg}{\sigma}            
\def\l{\ell}

\DeclareMathSymbol\crossrt{\mathrel}{AMSb}{"6E}
\def\suq2{U_q(su(2))}
\def\ot{\otimes}
\newcommand{\kett}[1]{|#1\rangle\!\rangle} 
\newcommand{\nn}{\nonumber}         
\newcommand{\CGq}[6]{C_q\!\begin{pmatrix}#1&#2&#3\\#4&#5&#6\end{pmatrix}}
\newcommand{\sesq}{{\mathchoice{\ooh}{\ooh}{\ssesq}{\ssesq}}}
\newcommand{\ssesq}{{\scriptstyle\frac{3}{2}}} 
\newcommand{\Dslash}{{D\mkern-11.5mu/\,}} 

\def\endproof{\unskip\nobreak\kern5pt\nobreak\vrule height4pt width4pt
depth0pt
\vskip4pt plus2pt}
\def\sq{\unskip\nobreak\kern5pt\nobreak\vrule height4pt width4pt
depth0pt}

\def\DO{{\CD}}

\def\CA{{\mathcal A}}
\def\CB{{\mathcal B}}
\def\CD{{\mathcal D}}
\def\CH{{\mathcal H}}
\def\CK{{\mathcal K}}

\def\CG{{\mathcal G}}
\def\cp{\Delta}

\def\ooh{\frac{3}{2}}

\def\ts{\otimes}

\def\acts{\triangleright}
\def\racts{\triangleleft}
\def\ket#1{\ee_{#1}}
\def\endproof{\vrule height 0.5em depth 0.2em width 0.5em}
\newbox\tbox
\newbox\aubox
\newbox\adbox
\def\title#1{\setbox\tbox=\hbox{\let\\=\cr
\baselineskip14pt\vbox{\Large\bf\tabskip 0pt plus15cc
\halign to\hsize{\hfil\ignorespaces \uppercase{##}\hfil\cr#1\cr}}}}
\newbox\abbox
\setbox\abbox=\vbox{\vglue18pt}
\def\author#1{\setbox\aubox=\hbox{\let\\=\cr
\baselineskip12pt\vbox{\tabskip 0pt plus15cc
\halign to\hsize{\hfil\ignorespaces \uppercase{{##}}\hfil\cr#1\cr}}}%
\global\setbox\abbox=\vbox{\unvbox\abbox\box\aubox\vskip8pt}}
\def\address#1{\setbox\adbox=\hbox{\let\\=\cr
\baselineskip12pt\vbox{\it\tabskip 0pt plus15cc
\halign to\hsize{\hfil\ignorespaces {##}\hfil\cr#1\cr}}}%
\global\setbox\abbox=\vbox{\unvbox\abbox\box\adbox\vskip16pt}}

\begin{document}
\thispagestyle{empty}
\begin{center}
{\LARGE
{\bf
GEOMETRY OF QUANTUM SPHERES
}
}
~\\
~\\~\\
{\Large
{\bf Ludwik D\c{a}browski}
~\\~\\
}
{\large
{\it
Scuola Internazionale Superiore di Studi Avanzati,\\
Via Beirut 2-4, I-34014, Trieste, Italy
}
\vspace{4cm}

{\bf
Abstract\\
}
}
\end{center}
{\large
\begin{quote}
Spectral triples on the q-deformed spheres of dimension two and three are reviewed.
\end{quote}
\begin{center}

\vspace{1cm}
\noindent {\it M.S.C.}: 81R60, 81R50, 20G42, 58B34, 58B32, 17B37. \\
~\\
\noindent {\it Key words and phrases}:\\
\begin{quote}
Noncommutative geometry, spectral triples, 
quantum groups, quantum spheres.\\
\end{quote}
\vspace{7cm}
Ref. SISSA 86/2004/FM\\
\end{center}
}

\newpage\noindent    
{\large
\section*{0. Introduction.}
The key notion of the most recent `layer' of noncommutative differential geometry 
a la` Connes \cite{C1} is {\em spectral triple}, which encodes the concept 
of a noncommutative Riemannian spin$_c$ manifold.\\
In the classical (commutative) situation the spectral triple $(\CA, \CH, D)$ 
canonically associated with a Riemannian spin$_c$ manifold 
consists of the algebra $\CA$ of smooth functions on $M$,
of the Hilbert ${\cal H}$ space of 
(square integrable) Dirac spinors, carrying a representation of $\CA$ 
(by pointwise multiplication), and of the Dirac operator $D$, 
constructed from the Levi-Civita` connection (metric preserving and torsion-free) 
plus a $U(1)$-connection.\\
If $M$ is spin and even-dimensional  
there exist the operators $J$ and $\gamma$ real structure and gradation
(known also as charge conjugation and parity). 
All these data are of great importance both in Mathematics and Physics.
They satisfy certain further seven properties which allow to reconstruct back
the underlying differential, metric and spin structure.\\
A generalization of these concepts to noncommutative algebras 
in the framework of noncommutative spectral geometry 
has already found plentiful applications. 
But the whole {\it zoo} of q-deformed spaces 
coming from the quantum group theory, 
was commonly believed not to match well the Connes' approach.
This was supported by some apparent ``no-go" hints such as that
exponentially growing spectrum of the quantum Casimir operator 
would prevent bounded commutators with the algebra,
some known differential calculi seemed not to come as bounded 
commutators with any $D$,
an early classification of equivariant representations
missed the spinorial ones
and also on some deformation theory grounds.
However, the intense recent activity
indicated a possibility to reconcile 
these two lines of mathematical research.\\
In this paper spectral triples on some of the simplest studied examples 
of quantum spheres of lowest dimension (2 and 3) are reviewed.
More precisely, we shall be concerned mainly with (the algebra of) 
the underlying space of the quantum group $SU_q(2)$ and its two homogeneous 
spaces known as the standard and the equatorial Podle\'s sphere
(see \cite{D} for the list of other low dimensional spheres).\\
~\\
For the sake of consistent as far as possible choice of conventions 
and notation for different examples some of the original formulae 
are given in an equivalent form.
In the sequel $0< q < 1$ and $ [x]=[x]_q$, where 
$ [x]_q :=\frac{q^x-q^{-x}}{q-q^{-1}}$ for any number $x$.

\section*{1. Spectral triples.}
We recall the general definition. \\
Definition 1. ~A (compact) ~{\bf spectral triple}  
~$(A, \CH, D)$~ consists of a (unital) {$*$\nobreak -\nobreak algebra} 
$A$ of bounded operators on a 
Hilbert space $\CH$ and a self-adjoint operator $D=D^\dagger$ on $\CH$ with\\ 
i) compact resolvent $(D-\lambda )^{-1}$, 
~ $\forall\lambda\in\Bbb C\setminus {\rm spec}D$\\
ii) bounded commutators ~ $[D, \a ]$,
~~$\forall \a\in A$.\\
A spectral triple 
is called {\bf even} if a grading $\gamma$ of 
$\CH$ is given,
$\gamma^2 = 1$, such that $\a\!\in\! A$ are even operators,
 $\a \gamma= \gamma \a$, and $D$ is odd $D \gamma   = - \gamma D$.
Otherwise it is called {\bf odd} (by convention then $\gamma=1$).\\
A spectral triple
is called {\bf real} if an anti-linear isometry
$J$ is given, 
whose adjoint action sends $A$ to its commutant
\begin{equation}
\label{tocomm}
[ \a, J\b J^{-1} ] = 0,
\quad \forall \a, \b \in 
A \ .
\end{equation}
We are particularly interested in spectral triples on homogeneous embeddable 
quantum spaces, described by a comodule subalgebra $A$ of a Hopf $*$-algebra $H$.
It is natural to employ their symmetry in order to reduce the 
search freedom or to find {\em equivariant triples}.
The formulation of the symmetry in terms of {\it coaction} and {\it coproduct} in $H$
has not yet been presented but one can work with the dual quantum group. \\
To deal with the equivariance we shall use {\it action} of a Hopf $*$-algebra 
$\widetilde U$.
Typically $\widetilde U$ is dual of $H$ in the sense of
nondegenerate Hopf algebra pairing
(e.g. a quantum universal enveloping algebra). 
This yields two commuting (left and right)
$\widetilde U$-module algebra structures on $H$
\begin{equation}
\label{actract}
u\acts h = h_1 ( u, h_2 ), ~~~~h \racts u = (u, h_1) h_2 \ ,
\end{equation}
where we use Sweedler's type notation for the coproduct in $H$.
The star structure is compatible with both actions:
\begin{equation}
u \acts h^* = ((Su)^* \acts h)^*,  \quad  h^* \racts u = (h \racts (Su)^*)^*,
\forall u \in U, \ h \in H,
\label{eq:star-act}
\end{equation}
where $S$ is the antipode in $\widetilde U$.\\
The actions (\ref{actract}) can be combined into a (left) action of 
$\widetilde U\ot \widetilde U^{op}$ via $(u\ot u^\prime ) ~{\acts}\  x := 
u ~{\acts}~ x~ {\racts}~ u^\prime $.
(Another option is to pass to the left action of $\widetilde U\ot \widetilde U$
by using $S$ and a suitable automorphism anti co-homomorphism of $U$).\\
In case $A$ is a proper subalgebra of $H$,
we assume that (at least) some nontrivial Hopf subalgebra $U$ of 
$\widetilde U\ot \widetilde U^{op}$ survives the restriction to $A$.\\
~\\
We shall use the following general definition.\\
Definition 2. ~Let $U$ be a Hopf $*$-algebra and $A$ be a $U$-module $*$-algebra.
A spectral triple $(A, \CH, D)$ is called {\bf $U$-equivariant} 
if there exists a dense subspace $V$ in $\CH$ such that\\
i) $V$ is a 
module over $A {\large \rtimes} U$, 
so in particular
\begin{equation}
\label{covar}
 u (\alpha v) = (u_{1} \acts \alpha ) (u_{2}~v)\
~~~\forall u \in U, \alpha\in A, v\in V\ ;
\end{equation}
ii) $u^*\subset u^\dagger$ ~on $V$ (as unbounded operators);\\
iii) $D(V)\subset V$ and $D$ is invariant
\begin{equation}
\label{Dinv}
Du = uD, ~~({\rm on} V), ~~\forall u \in U\ .
\end{equation}
Moreover, for $U$-equivariant even and/or real spectral triples
we shall require that $\gamma (V)\subset V$, $u \gamma = \gamma u$ and/or 
that $J$ is the phase part of some $\tilde J$, which satisfies
$\tilde J (V)\subset V$, $u \tilde J  = \tilde J (Su)^*$,  $\forall u \in U$.\\
~\\
(In (\ref{covar}) Sweedler's type notation is used for the coproduct in $U$
and $\acts$ for the action of $U$ on $A$).\\
~\\
On top of the above a series of further seven requirements (axioms) \cite{C2}
(see also \cite{G-BVF}) are required to describe a noncommutative 
(compact Riemannian spin) manifold, that we briefly recall.\\
A spectral triple $(A, H, D)$ is {\bf regular} (or {\bf smooth})
iff 
$$ [~|D|, \dots [~|D|, \b\underbrace{]\dots ]}_n$$
are bounded ~$\forall ~n\in{\Bbb N}, ~ \b\in {\cal A}\cup [D, {\cal A}]$.\\
This condition permits to introduce the analogue of Sobolev spaces
$\CH^s := Dom(|D|^s)$ for $s\geq 0$ 
(assume that ~$\CH^\infty := \cap_{s\geq 0} H^s$~ is a core of $|D|$).
Then $T:\CH^\infty \to \CH^\infty $ has {\em analytic order} $k$ iff 
$T$ extends to a bounded 
$T:\CH^{k+s}\!\to\!\CH^s$, $\forall s\geq 0$.\\
It turns out that 
~$\CA (\CH^\infty)\subset \CH^\infty$.
Moreover \cite{CM} (see also \cite{H})
certain algebra $\CD = \cup \CD_k$
of {\bf differential operators}  can be introduced
as the smallest algebra of operators on $\CH^\infty$
containing $A\cup [D, A]$ and 
filtered by the analytic order $k\in\Bbb N$ in such a way that
$[D^2, \DO_k] \subset \DO_{k+1}$.
Next the space $\Psi_k$ of {\bf pseudodifferential operators of order} 
$k\in\Bbb Z$ ~consist of 
those $T:\CH^\infty\to\CH^\infty$ which for any $m\in\Bbb Z$
~(especially $m<< 0$)~ can be written
in the form
$$T= S_0\,|D|^k + S_1\,|D|^{k-1} +\dots +  S_{k-m}\,|D|^m + R \ ,$$
where $S_j\in\DO_j $
and $R$ has analytic order $\leq m$.
(Assume that $|D|$ is invertible; otherwise work with $\sqrt{1+D^2}$).\\
It can be seen that
$\Psi_0$ is the algebra generated by the elements of $ A\cup [D, A]$
and their iterated commutators with $|D|$, and that
$[D^2,\Psi_{k} ]\subset \Psi_{k+1}$.
The algebra structure on 
$\Psi := \cup_{k\in\Bbb Z}\Psi_k$ 
can be read in terms of {\bf asymptotic expansion}:
$T \approx \sum_{j\in \Bbb N} T_j$,
whenever $T$ and $T_j$, for $j\in \Bbb N$, are operators
$\CH^\infty\to\CH^\infty$
and $\forall m\in \Bbb Z$,
~$\exists$ $N$ such that  $\forall M>N$, ~$T - \sum_{j=1}^M T_j$ ~
has analytic order $\leq m$.
For instance for complex powers of $|D|$ (defined by the Cauchy formula) one has
for $T\in\Psi$
$$ [~|D|^{2z}, T]  \approx
\sum_{j\geq 1}
{z\choose j}
[D^2, \dots [D^2, T\underbrace{]\dots ]}_{j}
|D|^{2z - 2j}\ .
$$
The algebra $\Psi$ provides a convenient framework to study
the residues of the zeta functions. For that it suffices
that {\bf dimension} requirement holds: $\exists$ $n\in \N$ s.t.
the eigenvalues (with multiplicity) of $|D|^{-n}$, 
$\mu_k = O(k^{-1})$ as $k\to \infty$. 
(The coefficient of the logarithmic divergence of 
$\sigma_N := \sum^N \mu_k$, denoted
$-\!\!\!\!\!\int  |D|^{-n}$, 
defines the `noncommutative itegral').\\ 
It follows that 
for $k > n$,
$D^{-k}$ is trace-class, i.e. {\bf finitely summable}.
Then ~{\bf dimension spectrum} is the set $\Sigma$
of the singularities of zeta functions
\begin{equation}
\zeta_\b (z) = {\rm Trace}_{\cal H} (\b\, |D|^{-z})
\end{equation}
for any $\b\in\Psi_0$.
~\\
If $\Sigma$ is assumed to be discrete with poles
only as singularities then
$\forall ~T\in\Psi_k$ the zeta function ${\rm Trace}(T|D|^{-z})$
is holomorphic in a half-plane of $\Bbb C$ with
$\Re z >>0$ and has a meromorphic continuation to $\Bbb C$
with all poles contained in $\Sigma + k$.
~\\
If in addition $\Sigma$ contains only simple poles
then the residue functional
$$\tau (T) := Res_{z=0} {\rm Trace}(T|D|^{-z})$$
is tracial on $T\in \Psi$ (c.f. \cite{Wo}).\\
These tools serve for the local index theorem of Connes-Moscovici \cite{CM},
which provides a powerful algorithm for performing complicated local computations 
by neglecting plethora of irrelevant details.\\
Another analytic requirement is {\bf finiteness and absolute continuity}:
$A$ admits a (Fr\'echet) completion ${\cal A}$ which is a pre $C^*$-algebra 
and $\CH^\infty$ is finite generated projective left ${\cal A}$-module.\\
~\\
Among the algebraic conditions for a spectral triple
we shall use the {\bf reality condition} 
iff  $(A, \CH, D)$ is real and certain sign conditions are satisfied 
for the square of $J$ and its (anti) commutation
with $\gamma$ and with $D$,
and the {\bf first order condition} iff  
\begin{equation}
\label{FOC}
\left[ [D, \a], J\b J^{-1} \right] =
0\;, ~~\forall \a, \b \in A\ .
\end{equation}
is satisfied.
There are two more algebraic conditions.
The first is {\bf orientability}: 
there exists a Hochschild $n$-cycle 
\begin{equation}
\sum c_0 \ot c_0^\prime \ot c_1 \ot ...\ot c_n) \in Z_n(\CA,\CA\ot\CA^{op}) \ ,
\end{equation}
such that 
\begin{equation}
\sum c_0J(c_0^\prime)^* J^{-1} [D, c_1]...[D, c_n] = \g \ ,
\end{equation}
where $n$ is the dimension and $\g$ the gradation operator.
The second one is {\bf Poincar\'e duality}, which can be formulated
as the requirement that the pairing 
\begin{equation}
K_0(\CA )\times K_0(\CA )\mapsto \Z, 
~\left([p], [q]\right) \mapsto 1/2\,  {\rm Index} 
\left(pJqJ^{-1} (1\! +\! \g D\g ) pJqJ^{-1} \right) 
\end{equation}
and \begin{equation}
K_1(\CA )\times K_1(\CA )\mapsto \Z, 
~\left([u], [v]\right) \mapsto 1/4\, {\rm Index} 
\left( (1\! +\! D/|D|) uJvJ^{-1} (1\! +\! D/|D|)\right) 
\end{equation}
is nondegenerate.

\section*{2. Quantum 3-sphere.}
In this section we present 
spectral triples on the quantum 3-sphere underlying the quantum group $SU_q(2)$.
We start by recalling its definition, 
the quantum symmetry algebra and its representations.\\
The unital $*$-algebra $A(SU_q(2))$ \cite{W} 
($0<q<1$) is generated by $\a, \b$ satisfying 
$$
\a \b = q \b \a \, , \ \a \b^* = q \b^* \a \, , \ \b \b^* = 
\b^* \b \ ,\\
\a^* \a + \b^* \b = 1 \, , \ \a \a^* + q^2 \b \b^* = 1 \, .
$$ 
The classical subset is (the `equator') $S^1$ given by the characters
$\b\mapsto 0, \a\mapsto \lambda$ with $|\lambda| = 1$.
The $C^*$-algebra $\CA$ associated to $A$ is isomorphic to the extension 
of $C(S^1)$  by $ {\cal K}\otimes C(S^1)$.
The $K$-groups are $K_0={\Z}$, $K_1=\Z$.\\
~\\
The dual (infinitesimal) symmetry of $A$ can be described using 
the quantized algebra $U_q(su(2))$, 
which is a $*$-Hopf algebra with generators $e,f,k,k^{-1}$ satisfying relations
\begin{equation}
ek = qke,  ~kf = qfk, ~k^2 - k^{-2} = (q-q^{-1})(fe-ef),
\nonumber
\end{equation}
and the coproduct
\begin{equation}
\cp k = k \ts k,  ~~~\cp e =  e \ts k + k^{-1} \ts e,
~~~ \cp f =  f \ts k + k^{-1} \ts f \ .
\nonumber
\end{equation}
Its counit $\epsilon$, antipode $S$, and $*$-structure are given
respectively by
\begin{equation}
\begin{array}{lllll}
 \epsilon(k) = 1,  &\phantom{xxx}&  \epsilon(e) = 0,
&\phantom{xxx}& \epsilon(f) = 0, \\
Sk = k^{-1}, &\phantom{xxx}& Sf = - qf, &\phantom{xxx}& Se = -q^{-1} e, \\
 k^* = k, &\phantom{xxx}& e^* =f, &\phantom{xxx}& f^* =e.
\end{array}
\nonumber
\end{equation}
The (commuting) left and right actions (\ref{actract}) of 
$U_q(su(2))$ on the generators of $A(SU_q(2))$ read explicitly
\begin{align}
k \lt \a &= q^{-\half} \a,    &  k \lt \a^* &= q^{\half} \a^*,  &
k \lt \b &= q^{\half} \b, &  k \lt \b^* &= q^{-\half} \b^*,
\nn \\
f \lt \a &= 0,            &  f \lt \a^* &= - \b^*,         &
f \lt \b &= q^{-1} a,            &  f \lt \b^* &= 0,
\label{Uact3} \\
e \lt \a &= q\b ,            &  e \lt \a^* &= 0,               &
e \lt \b &= 0,            &  e \lt \b^* &= - \a^*,
\nn
\end{align}
and 
\begin{align}
\a \rt k &= q^{-\half} \a,   &  \a^* \rt k &= q^{\half} \a^*,  &
\b \rt k &= q^{-\half} \b,   &  \b^* \rt k &= q^{\half} \b^*,
\nn \\
\a \rt f &= - b^*,      &  \a^* \rt f &= 0,               &
\b \rt f &= q^{-1}\a^*,         &  \b^* \rt f &= 0,
\label{Uract2}\\
\a \rt e &= 0,           &  \a^* \rt e &= q\b ,               &
\b \rt e &= 0,           &  \b^* \rt e &= - \a.
\nn
\end{align}
They can be combined into one (left) action of 
$U_q(su(2))\ot U_q(su(2))$ by 
$$(u\ot u') ~{\acts}\  x := 
u ~{\acts}~ x~ {\racts}~ (S^{-1} \Theta(u')) \ ,$$
where $\Theta: k\mapsto k^{-1}, e\mapsto -f, f\mapsto -e$.\\
We recall the irreducible finite dimensional representations
$\sg_\l$ of $U_q(su(2))$ labeled by {\em spin}
~ $\l = 0,\half,1,\sesq,2,\dots$ 
\begin{equation}
\begin{array}{l}
\label{repU}
 f\, \ket{\l,m} = [\l-m]^{\half}[\l+m+1]^{\half} ~\ket{\l,m+1}\\
 e\, \ket{\l,m} = [\l-m+1]^{\half}[\l+m]^{\half} ~\ket{\l,m-1}\\
 k\, \ket{\l,m}  = q^{m}\, \ket{\l,m}
\end{array}
\end{equation}
where the vectors 
$\ket{\l m}$, for $m = -\l, -\l\! +\!1,\dots, \l-1, \l$, 
form a basis for the irreducible $U$-module $V_\l$.
In fact $\sigma_\l$ are 
$*$-representations of $U_q(su(2))$, with respect to the hermitian
product on $V_\l$ for which the vectors $\ket{\l m}$ are
orthonormal.\\
~\\
We return to the quantum 3-sphere $SU_q(2)$.\\
Let $\chi$ be the normalized Haar state on $\CA$, 
and let $L^2(SU_q(2))$ be the Hilbert space associated with $\chi$.
Take the orthonormal basis of $L^2(SU_q(2))$ as
\begin{equation}
\label{basis}
\e_{\l,i,j}:=[2\ell +1]^{\frac{1}{2}} q^{-i} t^\l_{i,j}, ~
~\l \in 1/2\, {\N}, ~i,j \in \{ -\l , -\l\! +\! 1 , \ldots , \l \}\ ,
\end{equation}
where $t^\l_{ij}$ is $(i,j)$-th matrix element 
of the unitary irreducible corepresentation of spin $\l$.\\
\noindent
We shall use the left regular unitary representation of $\CA$ in $L^2(SU_q(2))$ 
which on the generators reads explicitly
\begin{equation}
\label{rep3a}
\a ~\e_{\l,i,j} = q^{2\l+i+j+1} \, \frac{(1-q^{2\l-2j+2})^{1/2} 
(1-q^{2\l-2i+2})^{1/2}}{ (1-q^{4\l+2})^{1/2} (1-q^{4\l+4})^{1/2}}
\, \e_{\l + \frac{1}{2}, i-\frac{1}{2} , j-\frac{1}{2}}\\
\end{equation}
$$
+ \frac{(1-q^{2\l+2j})^{1/2} 
(1-q^{2\l+2i})^{1/2}}{(1-q^{4\l})^{1/2} (1-q^{4\l+2})^{1/2}}
\, \e_{\l - \frac{1}{2}, i-\frac{1}{2} , j - \frac{1}{2}}
\nonumber
$$
\begin{equation}
\label{rep3b}
\b ~\e_{\l,i,j} = -q^{\l+j} \, \frac{(1-q^{2\l-2j+2})^{1/2} 
(1-q^{2\l+2i+2})^{1/2}}{(1-q^{4\l+2})^{1/2} (1-q^{4\l+4})^{1/2}}
\, \e_{\l + \frac{1}{2}, i+\frac{1}{2} , j-\frac{1}{2}}
\end{equation}
$$
+ q^{\l+i} \, \frac{(1-q^{2\l+2j})^{1/2} 
(1-q^{2\l-2i})^{1/2}}{(1-q^{4\l})^{1/2} (1-q^{4\l+2})^{1/2}} \, 
\, \e_{\l-\frac{1}{2}, i+\frac{1}{2} , j - \frac{1}{2}}\ .
$$
This representation satisfies {\em par excellence} the requirement 
ii) of Definition 2 with the full symmetry $U_q(su(2))\ot U_q(su(2))$.\\
~\\
We pass now to the question of spectral triples on $SU_q(2)$.\\ 
In \cite{BK} the Hilbert space of Dirac spinors is 
introduced as $\CH = {\IC}^2 \otimes L^2(SU_q(2))$ 
with the orthonormal basis chosen as
$$ v^{+}_{\l,ij}:= C^{\frac{1}{2},\l,\l+\frac{1}{2}}_{\frac{1}{2},j-\frac{1}{2},
 j} e_+ \otimes \e_{\l,i,j-\frac{1}{2}}+
 C^{\frac{1}{2},\l,\l+\frac{1}{2}}_{-\frac{1}{2},j+\frac{1}{2},
 j} e_- \otimes \e_{\l,i,j+\frac{1}{2}},$$
$$v^{-}_{\l,ij}:= C^{\frac{1}{2},\l,\l-\frac{1}{2}}_{\frac{1}{2},j-\frac{1}{2},
 j} e_+ \otimes \e_{\l,i,j-\frac{1}{2}}+
 C^{\frac{1}{2},\l,\l-\frac{1}{2}}_{-\frac{1}{2},j+\frac{1}{2},
 j} e_- \otimes \e_{\l,i,j+\frac{1}{2}};$$
where 
$\l \in \frac{1}{2} \, {\N}$, $i,j=-\l,-\l+1,...,\l; \l=0$, 
$C$ denote the $q$-Clebsch-Gordan coefficients \cite{BK}
and $e_+, e_-$ is the standard orthonormal basis in ${\IC}^2$
($e_\pm = \ket{\half, \pm\half}$ if identifing ${\IC}^2$ with $V_\half$).
The unitary representation of $\CA$ acts on $\CH$ as a tensor product of 
${\II}_2$ and the left regular representation (\ref{rep3a},\ref{rep3b}).
A candidate for the Dirac operator is defined in \cite{BK} by declaring
$v^{+}_{\l,i,j}, v^{-}_{\l,i,j}$ to be its eigenvectors with eigenvalues
$ [\l]_{q^2}$ and $ -[\l\! +\! 1]_{q^2}$ 
respectively. Unfortunately, it has unbounded commutators with the algebra $A$.\\
~\\
For that reason in \cite{G}
a modified operator $\tilde D$ has been proposed following the suggestion in 
\cite{CL} and \cite{CD-V}.
It has the same eigenvectors
$v^{+}_{\l,i,j}, v^{-}_{\l,i,j}$
but the corresponding eigenvalues 
$\l+\frac{1}{2}$ and $-(\l+\frac{1}{2})$ respectively,
have a linear growth.
It turns out that its absolute value $|\tilde D |$
~(though not $\tilde D$ itself) does satisfy the requirement of bounded commutators 
with the algebra $A$.
Thus, $(A , \CH , |\tilde D |) $  is a spectral triple.
Unfortunately it fails to capture topological 
information of the underlying noncommutative space
being $|\tilde D |$ a {\em positive} operator.
Indeed the sign of $| \tilde D |$ is trivial and the corresponding Fredholm module 
has trivial pairing with $K$-theory. 
Accordingly, the Poincar\'e duality axiom does not hold.\\
~\\
Another spectral triple on $SU_q(2)$ was presented in \cite{CP}.
Therein, the Hilbert space $\CH$ is just (one copy of) $L^2(SU_q(2))$ 
with orthonormal basis (\ref{basis}).
The unitary representation is just (\ref{rep3a},\ref{rep3b}).
The Dirac operator (up to a numerical factor) is 
$$
\label{D3}
D ~\e_{\l,i,j} = ( 1- 2 \, \d_{i,\l}) \, \l \,\, \e_{\l,i,j} \ ,
$$
where $\d$ is the Kronecker symbol. It has bounded commutators 
with $A$. Due to the term with $\d$ the condition iii) 
of definition 1 is satisfied only for one copy of $U_q(su(2))$, 
which is the `broken' symmetry of this spectral triple.
Though $D$ is not a positive operator, the Poincare` duality does not hold
\cite{CP3}.
Another feature is that there is no `classical' limit
since $[D,h]$ blows up for some $h\in A$ as $q\to 1$.
Moreover, $J$ is not provided and the first order condition for $D$
can not be established. In a sense, this example corresponds rather to a spin$_c$
manifold with Finsler metric.\\
Despite these peculiarities the spectral triple of \cite{CP}
has been intensively studied in \cite{C3} from the analytical point of view.
Therein, among other things, certain smooth subalgebra 
of $\CA$ was defined, stable under holomorphic functional calculus.
The regularity condition was verified and the pseudodifferential calculus
constructed. Its algebra of complete symbols was determined
(by computing the quotient by smoothing operators) and 
the cosphere bundle of $SU_q(2)$ constructed together with the geodesic flow.
The summability condition was verified and the dimension spectrum computed to be
$\Sigma = \{1, 2, 3\}$. The res$_0$ of the zeta function was explicitely computed
and the explicit  formula   for the local index cocycle obtained.\\
~\\
Yet another spectral triple on $SU_q(2)$ appeared in \cite{CP2}.
Therein, the Hilbert space $\CH$ has orthonormal basis $\e_{ij}$ with
$i\in\N$, $j\in\Z$. The (faithful) unitary representation is 
\begin{equation}
\label{rep3aa}
\a ~\e_{i,j} = (1-q^{2i})^{1/2}\, \e_{i-1,j}\ ,
\end{equation}
\begin{equation}
\label{rep3bb}
\b ~\e_{i,j} = q^{i}\, \e_{i,j-1}\ .
\end{equation}
It is equivariant under the action 
$\a\mapsto z\a, \b\mapsto w\b$ of the group $U(1)\times U(1)$,
implemented on $\CH$ by $~e_{i,j}\mapsto  z^{i}w^j\, e_{i,j}$.
A class of $U(1)\times U(1)$-invariant Dirac operators was identified
and among them
\begin{equation}
\label{DD}
D ~\e_{i,j} = \left( i\,{\rm sign}(j) + j\right)\, \e_{i,j}
\end{equation}
which is 2-summable and not positive.
In the `classical' limit $\a$ degenerates to 0 and $\b$ becomes the simple shift.  
This example also corresponds to a spin$_c$ manifold with Finsler metric
since $J$ is not provided and the first order condition for $D$
can not be established. However, interestingly the corresponding Connes-de Rham 
and the square integrable differential complexes have been computed.\\
~\\
Very recently a spectral triple appeared \cite{DLSSV}
with several nice properties.\\
The Hilbert space is $L^2(SU_q(2))\ot \IC^2$ 
(note the order !).
The representation of $U_q(su(2))\ot U_q(su(2))$
comes by coupling the first factor of 
$V = \bigoplus_{2l=0}^\infty  V_l \ox V_l$
with  $V_\half$.
So in particular, $u\ot u'$ is represented as
$$
\bigoplus_{2\l \geq 0} \sg_\l (u_1)\ot \sg_\l (u')\ot 
\sg_\half (u_2)\ .$$
Next, this coupling is diagonalized
by decomposing the tensor product with $\IC^2$ as 
$$\H = \H^\up \oplus \H^\dn = 
\bigoplus_{2j\geq 0} W_j^\up \oplus 
\bigoplus_{2j\geq 1} W_j^\dn \ .$$
(For consistency with the labeling of the representations by {\it spin} 
the index $\ell$ is offset by $\pm\half$).
Hence, the orthonormal basis is
\begin{equation}
\label{sbasisdn}
\ket{j\mu n\dn}
:= C_{j\mu} \,\ket{j^- \mu^+ n} \ox \ket{\half,-\half}
+ S_{j\mu} \,\ket{j^- \mu^- n} \ox \ket{\half,+\half};
\end{equation}
for the labels $j = l\! +\! \half$, $\mu = m\! -\! \half$, 
with $\mu = -j,\dots,j$ and $n = -j^-,\dots,j^-$ and
\begin{equation}
\label{sbasisup}
\ket{j\mu n\up}
:= - S_{j+1,\mu} \,\ket{j^+ \mu^+ n} \ox \ket{\half,-\half}
+ C_{j+1,\mu} \,\ket{j^+ \mu^- n} \ox \ket{\half,+\half},
\end{equation}
for the labels $j = l\! -\! \half$, $\mu = m\! -\! \half$, 
with $\mu = -j,\dots,j$ and $n = -j^+,\dots,j^+$.
Here the coefficients are
\begin{equation}
\label{sbasiscoeffs}
C_{j\mu} := q^{-(j+\mu)/2}\, [j-\mu]^\half [2j]^{-\half},  \qquad
S_{j\mu} := q^{(j-\mu)/2} \, [j+\mu]^\half [2j]^{-\half}
\end{equation}
and shorthand notation $ k^\pm := k \pm \half$ is adopted.
(Notice that there are no $\dn$ vectors for $j = 0$).\\
~\\
In order to simplify the formulae for the representation of $A$ on $\CH$,
the pair of spinors is juxtaposed as
\begin{equation}
\label{jsbasis}
\kett{j\mu n} := \begin{pmatrix} \ket{j\mu n\up} \\[2\jot]
\ket{j\mu n\dn} \end{pmatrix},
\end{equation}
for $j \in 1/2 \N$, with $\mu = -j,\dots,j$ and
$n = -j-\half,\dots,j+\half$, 
with the convention that the lower component is zero when
$n = \pm(j+\half)$ or $j = 0$. 
Furthermore, a matrix with scalar entries
$$
\tau = \begin{pmatrix} \tau_{\up\up} & \tau_{\up\dn} \\
\tau_{\dn\up} & \tau_{\dn\dn} \end{pmatrix},
$$
is understood to act on $\kett{j\mu n}$ by the rule:
\begin{align}
\tau \ket{j\mu n \up}
&= \tau_{\up\up} \ket{j\mu n\up} + \tau_{\dn\up} \ket{j\mu n\dn},
\nn \\
\tau \ket{j\mu n \dn}
&= \tau_{\dn\dn} \ket{j\mu n\dn} + \tau_{\up\dn} \ket{j\mu n\up}.
\end{align}
The unitary representation $\pi$ of $A$ on $L^2(SU_q(2))\ot \IC^2$ 
is the tensor product of 
the left regular representation (\ref{rep3a},\ref{rep3b}) and ${\II}_2$.
On the basis (\ref{jsbasis}) it reads explicitly
\begin{align}
a \,\kett{j\mu n}
&= \a^+_{j\mu n} \kett{j^+ \mu^+ n^+}
 + \a^-_{j\mu n} \kett{j^- \mu^+ n^+},
\nn \\[\jot]
b \,\kett{j\mu n}
&= \b^+_{j\mu n} \kett{j^+ \mu^+ n^-}
 + \b^-_{j\mu n} \kett{j^- \mu^+ n^-},
\nn \\[\jot]
a^* \,\kett{j\mu n}
&= \tilde\a^+_{j\mu n} \kett{j^+ \mu^- n^-}
 + \tilde\a^-_{j\mu n} \kett{j^- \mu^- n^-},
\label{eq:spin-repn} \\[\jot]
b^* \,\kett{j\mu n}
&= \tilde\b^+_{j\mu n} \kett{j^+ \mu^- n^+}
 + \tilde\b^-_{j\mu n} \kett{j^- \mu^- n^+},
\nn
\end{align}
where $\a^\pm_{j\mu n}$ and $\b^\pm_{j\mu n}$ are, up to phase factors
depending only on~$j$, the following triangular $2 \times 2$ matrices:
\noindent
\begin{align}
\a^+_{j\mu n} &= q^{(\half-\mu-n)/2}\, [j\! +\! \mu\! +\! 1]^\half
\begin{pmatrix}
q^{j+\half} \, \frac{[j+n+\sesq]^{1/2}}{[2j+2]} & 0 \\[2\jot]
q^{-\half} \,\frac{[j-n+\half]^{1/2}}{[2j+1]\,[2j+2]} &
q^{j} \, \frac{[j+n+\half]^{1/2}}{[2j+1]}
\end{pmatrix},
\nn \\[2\jot]
\a^-_{j\mu n} &= q^{(\half -\mu -n)/2}\, [j\! -\! \mu]^\half
\begin{pmatrix}
q^{-j-1} \, \frac{[j-n+\half]^{1/2}}{[2j+1]} &
- q^{-\half} \,\frac{[j+n+\half]^{1/2}}{[2j]\,[2j+1]} \\[2\jot]
0 & q^{-j-\half} \, \frac{[j-n-\half]^{1/2}}{[2j]}
\end{pmatrix},
\nn \\[2\jot]
\b^+_{j\mu n} &= q^{(-\mu-n-\half)/2}\, [j\! +\! \mu\! +\! 1]^\half
\begin{pmatrix}
\frac{[j-n+\sesq]^{1/2}}{[2j+2]} & 0 \\[2\jot]
- q^{j+1} \,\frac{[j+n+\half]^{1/2}}{[2j+1]\,[2j+2]} &
q^{\half} \, \frac{[j-n+\half]^{1/2}}{[2j+1]}
\end{pmatrix},
\label{eq:spin-coeff}
\\[2\jot]
\b^-_{j\mu n} &= q^{(-\mu -n-\half)/2}\, [j\! -\! \mu]^\half
\begin{pmatrix}
- q^{\half} \, \frac{[j+n+\half]^{1/2}}{[2j+1]} &
- q^{-j} \,\frac{[j-n+\half]^{1/2}}{[2j]\,[2j+1]} \\[2\jot]
0 & - \frac{[j+n-\half]^{1/2}}{[2j]}
\end{pmatrix},
\nn
\end{align}
and the remaining matrices are the hermitian conjugates
$$
\tilde\a^\pm_{j\mu n} = (\a^\mp_{j^\pm \mu^- n^-})^\7,  \qquad
\tilde\b^\pm_{j\mu n} = (\b^\mp_{j^\pm \mu^- n^+})^\7.
$$
This (spinorial) representation is equivariant 
(c.f. ii. of Definition 2) under the full symmetry 
$U_q(su(2))\ot U_q(su(2))$.\\
(Note that the decomposition of $\IC^2 \ot V$ as taken in \cite{BK} and \cite{G}
differs from $V \ox \IC^2$
by the $q$-Clebsch--Gordan coefficients according to the rule 
$$
\CGq j l m r s t = \CGq l j m {-s} {-r} {-t} ,
$$
which amounts to a substitution of $q$ by $q^{-1}$ in
the coefficients $C$ and $S$ (\ref{sbasiscoeffs}).
However, $\IC^2 \ot V$ is not equivariant in the sense of Definition 2.).\\
~\\
As far as the Dirac operator is concerned, it is diagonal with linear 
eigenvalues:\\
\begin{equation}
D\ket{j\mu n\up} = d_j^\up \,\ket{j\mu n\up},  \qquad
D\ket{j\mu n\dn} = d_j^\dn \,\ket{j\mu n\dn},
\label{eq:Dirac-eigen}
\end{equation}
where
\begin{equation}
d_j^\up = c_1^\up j + c_2^\up, \qquad
d_j^\dn = c_1^\dn j + c_2^\dn,
\label{eq:linear-evs}
\end{equation}
with $c_1^\up$, $c_2^\up$, $c_1^\dn$, $c_2^\dn$ independent
of ~$j$. (The multiplicities are $(2j + 1)(2j + 2)$ and $2j(2j + 1)$).\\
~\\
Such $D$ is selfadjoint on a natural domain and has a compact resolvent. 
The linearity of eigenvalues suffices to check that
$[D, x]$, $x\in A$, are bounded.
Moreover
$D$ is invariant under the full $U_q(su(2))\ot U_q(su(2))$ by construction.\\
In addition, if $c_1^\dn = - c_1^\up, ~c_2^\dn = - c_2^\up + c_1^\up$
then the spectrum of~$D$ coincides with that of the classical
Dirac operator $\Dslash$ on the round sphere $S^3$, 
up to rescaling and addition of a constant.
Thus, this spectral triple is an isospectral deformation of 
$(C^\infty(S^3),\H,\Dslash)$, and
in particular, its spectral dimension is~$3$.\\
Altogether, $(\CA(SU_q(2)),\CH, D)$ is a $3^+$-summable 
spectral triple, equivariant under the full symmetry $U_q(su(2))\ot U_q(su(2))$,
which is the rationale of its construction.\\
~\\
There is also an antilinear conjugation operator $J$ on~$\H$ which
is defined explicitly on the orthonormal spinor basis by
\begin{align}
J\, \ket{j\mu n\up} &:= i^{2(2j+\mu+n)} \,\ket{j,-\mu,-n,\up},
\nn \\
J\, \ket{j\mu n\dn} &:= i^{2(2j-\mu-n)} \,\ket{j,-\mu,-n,\dn}.
\label{eq:J-formula}
\end{align}
It is immediate that $J$ is antiunitary and
that $J^2 = -1$, since each $4j \pm 2(\mu + n)$ 
is an odd integer. 
Moreover, $J D J^{-1} = D$ since $D$ and $J$ diagonalize 
on their common eigenspaces $W_j^\up$ and $W_j^\dn$.\\
~\\
Besides these usual properties $J$ departs from the axiom scheme 
for real spectral triples.
Namely it does not map $\pi^\prime (A)$ to its
commutant, but it does (!) modulo the two-sided ideal $\CG$
of $\CB(\CH)$ generated by the (positive trace-class) operator 
$L_q \kett{j \mu n} := q^j \,\kett{j \mu n}$.
($\CG\subset \CK$ is contained in all ideals of infinitesimals of order 
$\alpha$ for any $\alpha >0$).\\
More precisely, for any $x,y \in A$,
\begin{equation}\label{3al-com}
\left[J x J^{-1}, y \right] \in \CG \, .
\end{equation}
Also the first order condition holds only up to $\CG$,
that is for all $x,y \in A$,
\begin{equation}
\left[ J x J^{-1}, [D, y] \right] \in \CG \ .
\label{3al-FOC}
\end{equation}
~\\
We remark that in a sense this spectral triple 
realizes the suggestion in \cite{CL} and \cite{CD-V}. 
However differently to \cite{G} (and \cite{BK}),
the spinor representation is constructed by tensoring 
$L^2(SU_q(2))$ by $\IC^2$ on the right rather than on the left.
Though the Clebsch--Gordan decomposition looks similar,
the different asymptotics of the appropriate Clebsch--Gordan coefficients 
leads to bounded commutator $[D,\pi (x)]$ instead of unbounded one 
for some $x$ in \cite{G} (where the off-diagonal terms in $\pi (x)$ 
are not compact but rather can be bounded from below). 
The origin of such a drastic difference 
is the (unbounded) non co-commutativity of $U_q(su(2))$.
\section*{3. Quantum 2-spheres.}
We pass now to 2-dimensional examples.
\subsection*{3.1. Standard quantum sphere.}
The algebra $A$ of the standard Podle\'s quantum sphere \cite{P1} can be defined 
in terms of generators $a \! =\! a^*, b, b^*$ and relations ($0\leq q <1$)
\begin{equation}
\begin{array}{lll}
\label{srel}
ab = q^2 ba, & \phantom{xxx} & ab^* = q^{-2} b^* a, \\
b b^* = q^{-2}a (1 - a), & & b^* b = a (1 - q^2 a).
\end{array}
\end{equation}
The associated $C^*$-algebra is isomorphic 
to the minimal unitization of the compacts $\cal K$.
It has one classical point given by the character
$a\mapsto 0, b\mapsto 0$ and $K$-groups are $K_0={\Z}^2$, $K_1=0$.\\
The algebra $A$ was discovered as a $A(SU_q(2))$-comodule algebra,
but it can be also viewed as a subalgebra of $A(SU_q(2))$ generated by
\begin{equation}
b=  -\a^* \b, ~~
b^* = -\b^* \a, ~~
a = \b\b^*=a^*.
\end{equation}
Thus the standard Podle\'s quantum sphere is not only a
homogeneous $SU_q(2)$-space, but also
a quotient space $SU_{q}(2)/{U(1)}$.
(In fact the quantum Hopf fibration and monopole connection are well known.)\\
~\\
The ~$\acts$~ action (\ref{Uact3}) of $U_q(su(2))$ descends to $A$,
up to equivalence it reads
$$
e \acts b= - (q^{\frac{1}{2}} + q^{-\frac{3}{2}}) a + q^{-\frac{3}{2}},
~e \acts b^* = 0,
~e \acts a = q^{-\frac{1}{2}} b^*,
\nonumber
$$
\begin{equation}
k \acts b = q b,
~k \acts b^* = q^{-1} b^*,
~k \acts a = a,
\label{acts}
\end{equation}
$$
f \acts b =0,
~f \acts b^* = (q^{\frac{3}{2}} + q^{-\frac{1}{2}}) a -
q^{-\frac{1}{2}},
~f \acts a = - q^{\frac{1}{2}} b.
\nonumber
$$
However this is not the case for ~$\racts$~. 
Indeed setting for $m\in \half\,\Z$
$${A}_m = \{h\in A(SU_q(2))~|~  h \, {\racts} \, k = q^{m}h \},$$
  (in particular ${A}_0 = A$)
  it is easy to see that ${A}_m{\racts} \, e = { A}_{m+1}$ ~and~
  ${ A}_m{\racts} \, f = { A}_{m-1}$.
  (Nevertheless ~$\racts$~  will still play a r\^ole in the sequel).\\
~\\
A spectral triple on the standard Podle\'s sphere has been constructed in \cite{DS}.
The Hilbert space $\CH$ has the orthonormal basis $ \ket{\ell,m,s}$ with
$\ell\in {\Bbb N} +1/2 $, ~$m\in\{-\ell, -\ell\! +\! 1, \dots , \ell\}$
~and ~$s\in\{-1, 1\}$.
The bounded representation $\pi$ of ${A}$ on $\CH$ reads:
\begin{equation}
\label{rep2}
\begin{array}{l}
a ~\ket{\ell,m,s} = a^+_{\ell,m,s}\ket{\ell\!+\! 1,m,s} +
a^0_{\ell,m,s}\ket{\ell,m,s} +  a^-_{\ell,m,s}\ket{\ell\!-\! 1,m,s},
\\
~\\
b ~\ket{\ell,m,s} = b^+_{\ell,m,s}\ket{\ell\!+\! 1,m\!+\! 1,s} +
b^0_{\ell,m,s}\ket{\ell,m\!+\! 1,s} +  b^-_{\ell,m,s}\ket{\ell\!-\! 1,m\!+\! 1,s},
\end{array}
\end{equation}
where
$$
a^+_{\ell,m,s}=
-q^{2\ell +m +1-s/2} \frac{(1\! -\! q^{2\ell-2m+2})^{1/2}(1\! -\! q^{2\ell+2m+2})^{1/2}}
{(1\! -\! q^{4\ell+4})^{1/2}(1\! +\! q^{2\ell+3}\! +\! q^{2\ell+1}
\! -\! q^{6\ell+7}\! -\! q^{6\ell+5}\! -\! q^{8\ell+8})^{1/2}}\ ,\hspace{10cm}
$$
$$
a^0_{\ell,m,s} =
(1+q^2)^{-1} + \hspace{15cm}
$$
$$
\quad\quad\quad
\frac{ 
 (1\! -\!  q^{4\ell+2}\!  +\!  q^{2\ell+2m}\!  +\!  q^{2\ell+2m+2})
\left( (1\! -\! q^{2\ell -1})(1\! -\! q^{2\ell+3})\!  -\!  sq^{2\ell +s}(1\! -\! q^2) 
\right)
}
{(1+q^2)(1- q^{4\ell})(1- q^{4\ell+4})}\ ,
$$
\begin{equation}
\label{rep2a}
a^-_{\ell,m,s}=
-q^{2\ell +m-1-s/2} \frac{(1- q^{2\ell-2m})^{1/2}(1- q^{2\ell+2m})^{1/2}}
{(1\! -\! q^{4\ell})^{1/2}(1\! +\! q^{2\ell+1}\! +\! q^{2\ell-1}
\! -\! q^{6\ell+1}\! -\! q^{6\ell-1}\! -\! q^{8\ell})^{1/2}}\ ,\hspace{10cm}
\end{equation}
$$
b^+_{\ell,m,s} =
q^{\ell -s/2} \frac{(1- q^{2\ell+2m+4})^{1/2}(1- q^{2\ell+2m+2})^{1/2}}
{(1- q^{4\ell+4})^{1/2}(1+ q^{2\ell+3}+q^{2\ell+1}
-q^{6\ell+7}-q^{6\ell+5}-q^{8\ell+8})^{1/2}},\hspace{10cm}
$$
$$
b^0_{\ell,m,s} =
 q^{m+\ell} \frac{(1\! - \! q^{2\ell+2m+2})^{1/2}(1\! -\! q^{2\ell-2m})^{1/2}
(sq^{2\ell +s}(1\! -\! q^2)\! -\! (1\! -\! q^{2\ell-1})(1\! -\! q^{2\ell+3}))}
{(1- q^{4\ell})(1- q^{4\ell+4})},\hspace{10cm}
$$
\begin{equation}
\label{rep2b}
b^-_{\ell,m,s} =
-q^{2m+3\ell -s/2} \frac{(1- q^{2\ell-2m})^{1/2}(1- q^{2\ell-2m-2})^{1/2}}
{(1- q^{4\ell})^{1/2}(1+ q^{2\ell+1}+q^{2\ell-1}
-q^{6\ell+1}-q^{6\ell-1}-q^{8\ell})^{1/2}}\ .\hspace{10cm}
\end{equation}
(The formulae for $b^*$ are similar).\\
%
The action of $U_q(su(2))$ on $A$ is as in (\ref{acts}).
The $*$-representation of $U_q(su(2))$ on $\CH$ is
a direct sum of the halfinteger spin ($\ell\in {\Bbb N} +1/2 $)
finite dimensional
irreducible representations $\sigma_\ell$ (\ref{repU}) of $U_q(su(2))$
with multiplicity 2 ($s=\pm 1$).
This corresponds to the simplest possible even spectral triple;~
higher multiplicities can be dealt accordingly.~\\
With this, $\pi $ is a direct sum of (the only) two inequivalent
$U_q(su(2))$-equivariant representations
$\pi_{\pm}$ of $A(S^2_q)$ on $\CH_\pm $ ($s=\pm 1$).
In fact, $\pi_{\pm}$ are irreducible representations of 
$A(S^2_q){\large \rtimes}U_q(su(2))$.\\
Next, in \cite{DS} the following operators are constructed:
the grading 
$$
\gamma \ket{\l,m,s} = s ~\ket{\l,m,s} \ ,
$$
the reality structure
\begin{equation}
J \ket{\l,m,s} = i^{2m}  \ket{\l,-m,-s}, \label{defJ}
\end{equation}
(antilinear isometry satisfying $J^2=-1$, $\gamma J = - J \gamma$
and sending $A(S^2_q)$ to its commutant),
and the Dirac operator
\begin{equation}
\label{dirac}
D \ket{\l,m,s} =  
[\l + 1/2 ]
~\ket{\l,m,-s},
\end{equation}
satisfying $D\gamma = -\gamma D$ and $D J = J D$ and 
$ D u =  u D,  \forall u \in U_q(su(2))$.
~\\
In \cite{DS} it is shown that such data 
$({\cal A}, {\cal H}, D, J, \gamma)$
define a (compact) real even $U_q(su(2))$-equivariant spectral triple,
which satisfies the first order condition.
Moreover, with the action (\ref{acts}) and the representation 
(\ref{rep2},\ref{rep2a},\ref{rep2b}), 
$D$ is unique up to multiplication of the r.h.s of (\ref{dirac}) 
respectively by $z$ (or $\bar z$) when $s= +1$ (or $s= -1$), 
where $z\in {\Bbb C}\setminus \{ 0\}$.
Evidently $D$ is a (unbounded) selfadjoint operator defined 
on a natural dense domain in $\CH$ 
consisting of $\psi =\sum_{\l, m, s} c_{\l, m, s}  \ket{\l,m,s}~$,
$c_{\l, m, s}\in {\Bbb C}$,
such that
$\sum_{\l, m, s} \left(1+ 
[\l + \frac{1}{2}]^2 \right)
|c_{\l, m, s}|^2  < \infty$.\\
The limit $q\to 1$ is in agreement with the classical 
spectral triple on $S^2$,
identifying $ \ket{\l,m,_\pm}$ with $Y^\pm_{\l,m}$.
Also,  $J$ coincides with the charge conjugation on spinors
and $D$ has the correct $q=1$ limit.\\
The geometrical meaning of $\CH$, $J$, $D$ given above is quite clear.
Recall that classically Dirac spinor $\psi$ on $S^2$ is a
section of certain rank two vector bundle
associated to the bundle of {\em spin frames}, 
which is nothing but the Hopf bundle $U(1)\to SU(2)\to S^2$.
Equivalently, $\psi$ can be thought of as $U(1)$-equivariant 
(under $z \to z^{-1} \oplus z$) ${\IC}^2$-valued function on $SU(2)$, 
~or equivalently, as $-\half \oplus \half$ eigenvector of ~$\racts H$ 
(Cartan generator of $su(2)$).
Of course the space $\Gamma = \{\psi \}$ of such polynomial spinors 
is a finite projective module over $A$.\\
~\\
The  q-deformation of these bundles (and more) is well known \cite{BM},\cite{HM} 
and we can view $SU_q(2)$ as the bundle of spin q-frames, and Dirac spinors 
$\Gamma_q$ as $q^{-\half} \oplus q^{\half}$ eigenspace of 
the action ~$\racts $~ of $K$ on $A(SU_q(2))$.
Explicitly, 
$$\Gamma_q = A_{-\half} \oplus A_{\half}=
{\rm span}_A \{\a, \b \} \oplus {\rm span}_A \{\a^*, \b^*\}\ .$$
It is a finite projective module over $A$, with a projector in
${\rm Mat}(2, A)\oplus {\rm Mat}(2, A)$ being
%
\begin{equation}
\left(
\begin{array}{ll}
\a \a^*&\a \b^*\\ \b \a^*&\b \b^*
\end{array}
\right)
\oplus
\left(
\begin{array}{ll}
\a^*\a &q\a^*\b\\ q\b^*\a &q^2\b^*\b
\end{array}
\right)
\end{equation}
\begin{equation}
=\left(
\begin{array}{ll}
1\! -\! q^2a&~-qb^*\\-qb&~a
\end{array}
\right)
\oplus
\left(
\begin{array}{ll}
1\! -\! a&~-qb\\ -qb^* &~q^2 a
\end{array}
\right).
\end{equation}
The Hilbert scalar product comes by restriction from $L^2 (SU_q(2))$
and sending $ \ket{\ell,m,s} \to [2\l +1]^{- 1/2} q^{-s/2}\e_{\ell,s/2,m}$ 
gives unitary isomorphism.
This explains the relation of $\CH_\pm$ and $\CH$ to q-deformed Hopf bundles.\\
With this setup the above Dirac operator is just
$\left(
\begin{array}{cc}
0&\racts \, e  \\\racts f&0
\end{array}
\right)
$.
Then the invariance (\ref{Dinv}) is clear and 
since $\racts e$ and $\racts f$ are twisted derivations
it is immediate to see \cite{SW}, \cite{NT} that for $x\in A$ 
\begin{equation}
\label{[Dh]}
[D, x] =
\left(
\begin{array}{cc}
0&q^{1/2}h\racts e   \\q^{-1/2}h\racts f&0
\end{array}
\right)
\end{equation}
are bounded.
~\\
We comment on other axioms.
Concerning the dimension axiom the eigenvalues
of $|D|$ are $\pm [k]_q$ with degeneracy $4k$, 
where $k := \ell\! +\!  1/2 \in \Bbb N$.
Thus the eigenvalues of $|D|^{-z}$
decrease exponentially as $k\to \infty$ ~for $Re z\! >\! 0$ (recall that $q<1$).
By summing the derivative of the geometric series
finite-summability holds, actually $\epsilon$-summability for any $\epsilon >0$
$$
{\rm Trace}_{\cal H} |D|^{-\epsilon} <
4 (q-q^{-1})^{\epsilon}q^\epsilon/ (q^\epsilon -1)^2\ .
$$
Next $\sigma_N / log N \to 0$ for all $z > 0$,
$\to \infty$ for $z \leq 0$.
So, loosely speaking, the metric dimension of this geometry is `$0_+$',
which matches the known drop by 2 of the cohomological dimension of q-spaces.
More precisely, the dimension spectrum $\Sigma$ should be studied,
e.g. the singularities of
$
\zeta_\b (z) = {\rm Trace}_{\cal H} (\b\, |D|^{-z})
$.
For $\b =1$ its most divergent part goes as
$$
\sim 4 (q-q^{-1})^z \sum_{k\in\Bbb N } k q^{z k}
~\sim~
~ \frac{4 (q-q^{-1})^{2\pi in/h}}{ (z-2\pi in/h)^{2}}~
$$
as $z\sim 2\pi in/h$, where $h=\log q$, $\forall n\in\Bbb N$ 
(at least for Re(z)$\geq 0$). 
This suggests that $\Sigma$ should 
contain a (discrete) infinite lattice ~$2\pi in/h$,
$n\in\Bbb Z$ of {\em double} poles on the imaginary axis.\\
That simple poles are excluded could be expected
from the divergence of $\sigma_N / log N$ at $z=0$,
due to Tauberian theorem. Recall also that 
complex poles may appear on fractals and multiple poles on singular manifolds.\\
~\\
Concerning the reality requirement $J$ is tailored as in the (naive) dimension $d=2$,
(as indicated by the signs in $J^2 = -1, JD=+ DJ$ and $J\gamma = - \gamma J$),
which has at least the same parity as metric dimension.\\
~\\
The regularity axiom does not hold \cite{NT} in its rigorous form
but it looks as if there might be some substitute for it 
since for $h\in A_m$, ~$|D|^z h |D|^{-z} - q^{mz}h$
is of order -1.
Thus, by (\ref{[Dh]}) and recalling that $h\racts e\in A_1$ and
$h\racts f\in A_{-1}$ for $h\in A_0 =A$, it appears 
that the principal symbol of $|D|$ is not a scalar.
But one can still work with this sort of `spin anomaly' by using q-mutators since
$$|D|^{-z}[D,h ] -
{\small\left( \begin{array}{ll} q^z & 0 \\ 0 & q^{-z} \end{array} \right)}
[D,h] |D|^{-z} = \kappa (z) |D|^{-z-1}\ ,$$
where $\kappa (z)$ are bounded operators analytic in $z$.
In fact a `twisted local index formula' for $SU_q(2)$
was worked out in \cite{NT} in the framework of twisted cyclic cohomology.\\
~\\
Concerning the orientation axiom
perhaps some modification of \cite{SW}
containing some proofs due to Heckenberger, and of \cite{NT}
may be helpful.
They show, that $ d \a := i[D, \pi(\a)]$ gives
the unique 2 dimensional covariant first order differential calculus of \cite{P2}.\\
The related universal differential calculus contains a $A$-central 2-form,
with which one can associate (via the Haar measure $\chi$, such that 
$\chi(ab)=\chi(b\sigma(a))$, where $\sigma$ is the algebra automorphism of 
$A$ given by $\sigma  = K^2 \acts\, $)
a non-trivial {\em $\sigma$-twisted} cyclic Hochschild 2-cocycle $\tau$
on $A$
$$
\tau (\a,\b,\g) :=
\chi\left( \a(\b\racts\, e)(\g\racts f) - q^2\a(\b\racts f)(\g\racts\, e)\right)$$
$$
= {(q\! -\! q^{-1})^{-1}}{\log q}~
{\rm Res}_{z=2}~ {\rm Tr}_{\cal H}\,
{\small\left( \begin{array}{ll} 1 & 0 \\ 0 & -q^2 \end{array} \right)}
 K^2 |D|^{-z} \a [D, \b][D, \g]
$$
\cite{SW} (see also versions 1,2 of \cite{G}).
\subsection*{3.2. Equatorial quantum sphere.}
The algebra of equatorial Podle\'s sphere can be written as
$$ 
ba = q^2 ab, ~~b^*b + a^2 =1, ~~bb^* + q^4 a^2 = 1.
$$
The classical subset is (the `equator') $S^1$ given by the characters
$a\mapsto 0, b\mapsto \lambda$ with $|\lambda| = 1$.
The related $C^*$-algebra is extension of $C(S^1)$ by $ \cal K\oplus \cal K$
and $K_0={\Z}^2$, $K_1=0$.\\
The symmetry used for the equivariance is now
$$k \acts a =  a,  ~e \acts a =  q^{-\frac{1}{2}} {b^*}, 
~f \acts a = -q^{-\frac{3}{2}} b,
$$
$$
k \acts b = q^{-1} b, ~e \acts b = - (1+q^{-2}) q^{\frac{5}{2}} a, ~f \acts b =  0,
$$
$$
k \acts b^* = q {b^*}, ~e \acts {b^*} = 0, 
~f \acts {b^*} =  (1+q^{-2}) q^{\frac{3}{2}} a.
$$
In \cite{DLPS} the data $\CH $, $\gamma$ and $J$ 
are as in sect. 3.1 and the representation (up to equivalence) is
$\pi = \pi_+\oplus \pi_-$ where $\pi_\pm$ are
the (only) two irreducible $\suq2$-equivariant $*$-representations
$$
a~ \ket{\ell,m} = 
- q^{\ell-m}\frac{(1-q^{2\ell-2m+2})^{1/2}(1-q^{2\ell+2m+2})^{1/2}}{1-q^{4\ell+4}} 
\;\ket{\ell\! +\! 1,m} 
$$
\begin{equation}
\label{rep2ae}
\mp q^{2\ell-1}\frac{(1-q^2)(1+q^{4\ell+2}-q^{2\ell-2m}-q^{2\ell-2m+2})}
{(1-q^{4l})(1-q^{4l+4})}
\; \ket{\ell,m} 
\end{equation}
$$
- q^{\ell-m-1}
\frac{(1-q^{2\ell-2m})^{1/2}(1-q^{2\ell+2m})^{1/2}}{1-q^{4\ell}} \; \ket{\ell\! -\! 1,m},
$$
$$ 
b~ \ket{\ell,m}  = q^{2\ell-2m+1}
\frac{(1-q^{2\ell+2m+2})^{1/2}(1-q^{2\ell+2m+4})^{1/2}}{1-q^{4\ell+4}} \;
\ket{\ell\! +\! 1,m\! +\! 1}
$$
\begin{equation}
\label{rep2be}
 \pm 
q^{3\ell-m-1}\frac{ (1-q^4)(1-q^{2\ell+2m+2})^{1/2}(1-q^{2\ell-2m})^{1/2}}
{(1-q^{4l})(1-q^{4l+4})} 
\; \ket{\ell,m\! +\! 1} 
\end{equation}
$$
- \frac{(1-q^{2\ell-2m})^{1/2}(1-q^{2\ell-2m-2})^{1/2}}{1-q^{4\ell}} \;
\ket{\ell\! -\! 1,m\! +\! 1}.
$$
They are equivalent to the representation obtained by restricting 
the representation (\ref{rep3a},\ref{rep3b}) of $A(SU_q(2))$ to the subalgebra $A$
and by restricting $\cal H$ to the ($L^2$-completion) 
of certain vector spaces (left $A$-modules) constructed in \cite{BM2}.\\
~\\
It is shown in \cite{DLPS} that $J$ does not map $\pi(A)$ to its
commutant, but it does modulo the ideal $\CG\subset \CK$
generated by the operator $q^l$ on $\ket{\l,m,s}$.
More precisely, for any $x,y \in A$,
\begin{equation}\label{al-com}
\left[J xJ^{-1}, y \right] \in \CG \, .
\end{equation}
The operator $D$ is 
\begin{equation}
\label{dirac2}
D \ket{\l,m,s} =  (\l + 1/2 ) ~\ket{\l,m,-s} \ .
\end{equation}
It is unique (up to rescaling and addition of a constant) 
under the same postulates as in sect. 3.1 
except the first order condition is required only up to $\CG$,
that is for all $x,y \in A$,
\begin{equation}
\left[ J x J^{-1}, [D, y] \right] \in \CG \ .
\label{al-FOC}
\end{equation}
\noindent
More importantly, for any $x \in A$, the commutators $[D, x]$ are bounded.
Moreover, it is evident that $D$ is self-adjoint on a natural
domain in $\CH$ and has compact resolvent. 
Since the eigenvalues of $|D|$ are
$k=\l+\frac{1}{2} \in{\N}$ with multiplicity $4k$ 
the deformation is isospectral and the dimension requirement is satisfied with
the spectral dimension of ($A, \CH, D)$ being $n=2$.

\subsection*{3.3. Generic quantum sphere.}
The algebra of generic Podle\'s sphere $S_{q,c}^2$ can be written as
$$ 
ba = q^2 ab, ~~b^*b + a^2 -a = c1, ~~bb^* + q^{4} a^2 - q^{2} a = c1\ ,
$$
where $q<1$, $0<c<\infty$.
The classical subset is $S^1$ given by the characters
$a\mapsto 0, b\mapsto \lambda$ with $|\lambda| = c$.
The related $C^*$-algebra is extension of $C(S^1)$ by $ \cal K\oplus \cal K$
and $K_0={\Z}^2$, $K_1=0$.\\
~\\
A spectral triple on $S_{q,c}^2$ appeared in \cite{CP2}.
Therein, the Hilbert space $\CH$ has orthonormal basis $e_{n,s}$ with
$n\in\N$, $s\in\{-1,+1\}$. The (faithful) unitary representation is 
\begin{equation}
\label{rep3aaa}
a ~\ee_{n,s} = \left( 1/2\!  +\!  s(c\! +\! 1/4)^{1/2}\right) q^{2n}\, \ee_{n,s}\ ,
\end{equation}
\begin{equation}
\label{rep3bbb}
b ~\ee_{n,s} = \left( (1/2 + s(c\! +\! 1/4)^{1/2}) q^{2n}    
\! -\! (1/2\! +\! s(c\! +\! 1/4)^{1/2})^2 q^{4n}\! +\! c\right)^{\frac{1}{2}}\ee_{n,s}.
\end{equation}
It is equivariant under the action 
$\a\mapsto z\a, \b\mapsto w\b$ of the group $U(1)\times U(1)$,
implemented on $\CH$ by $~e_{i,j}\mapsto  z^{i}w^j\, e_{i,j}$.
A class of $U(1)\times U(1)$-invariant Dirac operators was identified
and among them,
\begin{equation}
\label{DDD}
D ~\ee_{n,s} = n\, \ee_{n,-s}\ ,
\end{equation}
which is 1-summable, not positive and has bounded commutators with the algebra.
This spectral triple is even, $\gamma \ket{n,s} = s ~\ket{n,s}$, 
has nontrivial Chern character but $J$ and the first order condition 
are not considered and thus it corresponds to a spin$_c$ manifold with Finsler metric. 
However, interestingly the corresponding Connes-de Rham differential complex 
has been computed.

\section*{4. Final comments.}
It is worth mentioning that $D$, or its spectrum, taken individually 
carries only part of the available geometric information.
It is the interplay of all the data of the spectral triple
(also $J$ and $\gamma$)
on the Hilbert space, that imposes some really stringent
restrictions and produces the spectacular consequences.\\
~\\
The examples of some recently constructed 
spectral triples on quantum spheres we have described
dissolve the widespread belief that Connes' approach to noncommutative geometry
does not match quantum-group theory.\\
The last example in Section 2 acting on two copies of the $L_2$ space
is odd, isospectral to the classical Dirac operator, 3-summable
and has the $U_q(sl(2))\otimes U_q(sl(2))$ symmetry.
It exhibits a variation of `up to infinitesimals' 
of the reality condition and the first order condition,
which was anticipated by the even, isospectral, 2-summable, 
$U_q(sl(2))$-equivariant example in section 3.2.
The (even) example of sect. 3.1 is not isospectral but it satisfies the 
reality and the first order axioms in the strict sense.\\
We believe the isospectral deformations will prove useful in the $q$-geometry
and will be omnipresent on other $q$-deformed spaces.
In fact this point is is currently studied on general Podle\'s sphere 
(including the standard one).\\
~\\\noindent
There are several interesting open mathematical problems
regarding the analytic properties of these spectral triples
(regularity and finiteness)
as well as the the algebraic ones (orientation and Poincar\'e duality).
The whole analysis along the lines of \cite{C3} 
culminating with the local index formulae is in preparation.\\
From the more physical point of view the interesting questions regard the 
construction of a wider class of spectral triples 
and of the Yang-Mills and gravitation theories on q-deformed spaces.
Does the action functional of \cite{CC}, \cite{CC2} attain the extrema
on the most symmetric spectral triples or is sensitive  
to the value of $q$ or of other parameters?

\section*{Acknowledgments}
This work was partially supported by the  EU Project INTAS 00-257.


}
\end{document}